\documentclass{article}

\usepackage{arxiv}

\usepackage[utf8]{inputenc} 
\usepackage[T1]{fontenc}    
\usepackage{hyperref}       
\usepackage{url}            
\usepackage{booktabs}       
\usepackage{amsfonts}       
\usepackage{nicefrac}       
\usepackage{microtype}      
\usepackage{graphicx}
\usepackage{amsmath}
\usepackage{amsthm,bm}

\graphicspath{ {./images/} }

\usepackage{threeparttable}  
\usepackage{enumitem}

\title{Fast and Accurate Reconstruction of Voronoi Generators in Large Tessellations}

\author{
 Carlos M Hernandez-Suarez \\
  CGIC\\
  Universidad de Colima\\
  Colima, Colima 15213 \\
  M\'exico \\
  \texttt{cmh1@cornell.edu} \\}
   
\begin{document}
\maketitle

\begin{abstract}
A Voronoi diagram partitions the plane into convex cells, each containing the points closest to a single generator.  
Given such a tessellation, the inverse Voronoi problem seeks the generator set \(S\) that produced it.  
Our algorithm selects a single interior cell with \(k\) edges and solves a compact, consistent linear system with \(2(k+1)\) unknowns and \(4k\) scalar equations to recover that cell’s generator together with the \(k\) generators of its neighbors in one step; the remaining sites follow by successive geometric reflections.  
The overall running time is \(\mathcal{O}(n)\) for a diagram with \(n\) cells.  
Across \(10^{3}\) Monte Carlo simulations on diagrams of \(10^{4}\) cells, the method achieved an average RMSE of \(10^{-12}\) and a worst-case individual reconstruction error of \(10^{-8}\), demonstrating both efficiency and robustness.
\end{abstract}


\section{Introduction}

Voronoi diagrams play a fundamental role in a wide range of scientific and engineering disciplines, including spatial ecology, wireless communications, urban planning, and computational geometry.  
Given a finite set of generator points (or \emph{sites}) in the Euclidean plane, a Voronoi diagram partitions the plane into convex cells, each consisting of all points closer to one site than to any other.

The \emph{Voronoi inverse problem} asks the converse question:  
given only the tessellation, recover the underlying set of generator points.  
This situation occurs, for instance, when a cellular coverage map or a biological mosaic is observable, but the locations that induced the partition are not.

Central to every reconstruction technique is the following elementary yet powerful geometric fact. For a full proof and further background, see Aurenhammer and Klein~\cite{aurenhammer2000voronoi}:

\medskip
\noindent\textbf{Lemma 2.1 (Perpendicular-bisector property).}  
\emph{For any two adjacent cells whose sites are \(P_i\) and \(P_j\), the edge common to both cells lies on the perpendicular bisector of the segment \(P_iP_j\).}

Lemma 2.1 converts geometric adjacency information into linear constraints on the (unknown) site coordinates and therefore underpins the algebraic formulations surveyed below.

First, Evans and Jones \cite{evans1987detecting} employed Lemma 2.1 to derive a linear system with \(2K\) equations in \(2n\) unknowns, where \(n\) is the number of sites and \(K\) the number of edges in the tessellation.  
Subsequently, Hartvigsen \cite{hartvigsen1992recognizing} demonstrated that the resulting program can be solved in polynomial time using linear programming.

Iterative geometric approaches constitute a second family of methods.  
Adamatzky \cite{adamatzky1993massively} proposed a massively parallel algorithm that places a provisional generator in a chosen cell and recursively reflects it across cell edges to reveal neighboring sites, continuing outward until the entire diagram is reconstructed. For arbitrary planar straight-line graphs, a linear-size fitting algorithm is given by Aloupis et al. \cite{Aloupis2013}; their method inserts \(\mathcal{O}(n)\) new sites to ensure every original edge appears as a Voronoi edge.  Unlike our reconstruction, which recovers the unique generator of each cell when the diagram is genuinely Voronoi, their algorithm trades accuracy for universality and site minimization.

\begin{figure}[ht]
  \centering
  \includegraphics[width=0.5\textwidth]{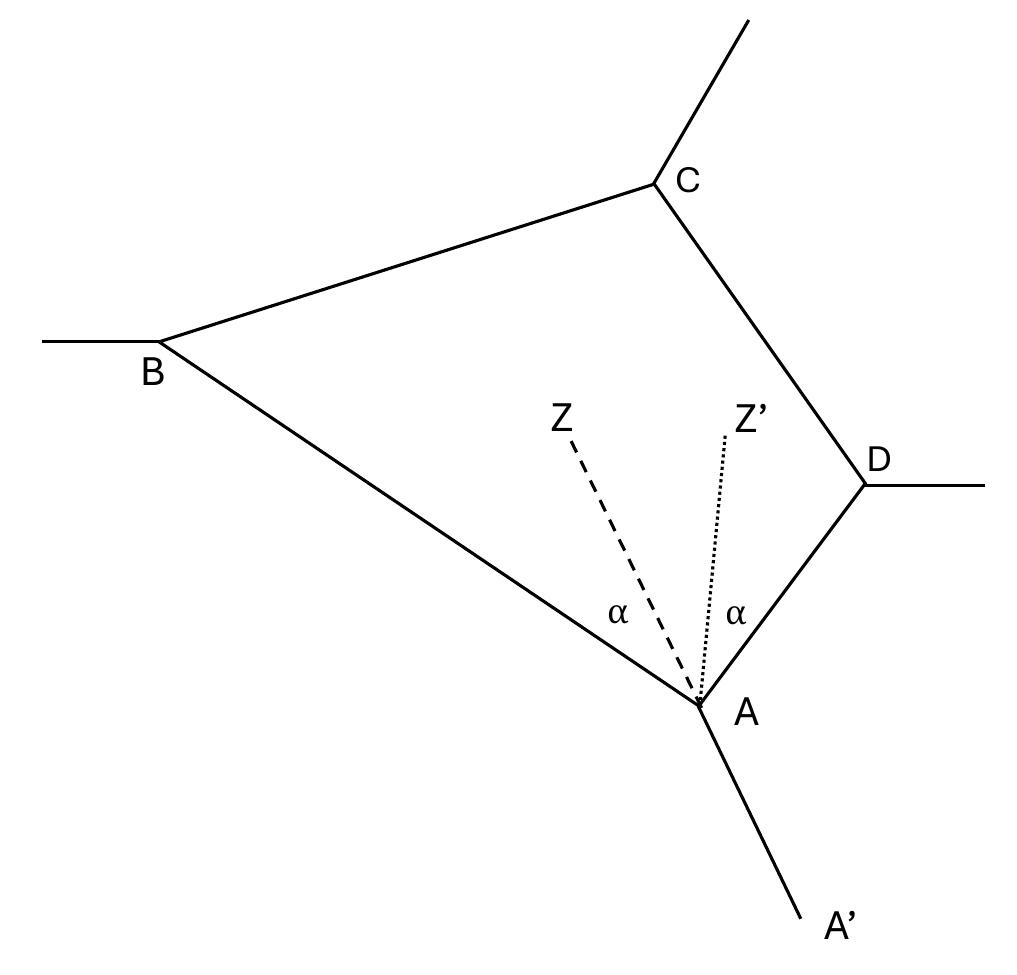}
\caption{Illustration of Schoenberg's \cite{schoenberg2003inverting} method:  $AZ$ is the extended outer ridge $AA’$ into the cell. This ray meets the nearest side $AB$ at angle $\alpha$.  Rotating $AZ$ by $\alpha$ gives $AZ’$, which passes through the generator.  Repeating the step with another outer ridge produces a second ray; their intersection pinpoints the generator.}
  \label{fig:fig1}
\end{figure}

An influential line of work is due to Schoenberg, Ferguson, and Li \cite{schoenberg2003inverting}.  
Their construction extends a boundary ridge into the interior of the cell, rotates the ray by the incident angle, and repeats the procedure with a second ridge; the intersection of the two rotated rays yields an estimate of the site (see Figure \ref{fig:fig1}).  
Six variants were examined, and the most accurate—denoted \(C'\)—has become a de facto benchmark.  
Notably, \cite{schoenberg2003inverting} is among the very few studies to report systematic Monte Carlo simulations, providing quantitative reference data against which new algorithms, including the present one, can be compared.

Our method can be summarised as follows.  
Starting from a single strictly interior \emph{anchor} cell with \(k\) edges, we assemble one compact, over-determined but consistent linear system whose \(2(k+1)\) unknowns comprise the coordinates of that cell’s site together with those of its \(k\) neighbors.  
Solving this system once yields those first \(k+1\) sites simultaneously; the remaining sites follow by rigid reflections.  
The resulting algorithm runs in \(\mathcal{O}(n)\) time for a diagram with \(n\) cells, attains machine-precision accuracy in large-scale Monte Carlo tests, and outperforms the benchmark \(C'\) method of \cite{schoenberg2003inverting} by one to two orders of magnitude in reconstruction error.

\section{Methods}
\label{sec:methods} 
\begin{figure}[ht] \centering \includegraphics[width=0.5\textwidth]{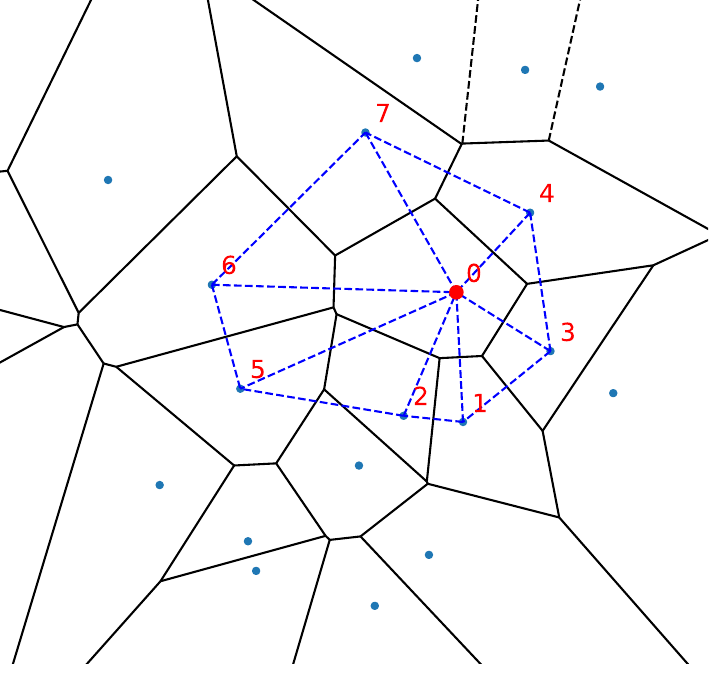} \caption{All generators of cells adjacent to the anchor cell~0 are reflections of its generator \(\mathbf g_{0}\). Moreover, every neighboring pair \(j,k\) in this ring are reflections of one another. These constraints give rise to a linear system with \(2(k+1)\) unknowns and \(4k\) scalar equations.} \label{fig:fig2} \end{figure}

Traditional inverse-Voronoi algorithms proceed in two stages. First, the generator of a single cell is recovered; second, that local procedure is either repeated for every remaining cell or the recovered site is successively reflected into adjacent cells. Our approach differs fundamentally: we select a non-degenerate, strictly interior cell with \(k\) edges and assemble a compact linear system whose \(2(k+1)\) unknowns are the coordinates of that cell’s generator together with those of its \(k\) neighbors. The system is over-determined (\(4k\) scalar equations versus \(2(k+1)\) unknowns) yet consistent and uniquely solvable. Solving it once yields the entire \((k+1)\)-site patch in a single computation, thereby avoiding repeated local reconstructions or iterative refinements. The generators of all remaining cells are then obtained successively by rigid reflections. Lemma 2.1 allows all generators in a Voronoi tessellation to be derived by geometric reflection. Let \(T_{ij}\) denote the line containing the ridge separating cells \(i\) and \(j\). Fix a generator \(\mathbf g_{0}\); then the generators of all adjacent cells are the mirror reflections of \(\mathbf g_{0}\) across the corresponding lines \(T_{0j}\). Moreover, for any two adjacent cells \(j\) and \(k\) in this neighborhood, their generators are mirror reflections of one another across the line \(T_{jk}\) (see Figure \ref{fig:fig2}).

\subsection{First local reconstruction around an anchor cell} Let \(A\) be a finite, non-degenerate \emph{anchor cell} and \(B=\{B_{1},\dots,B_{k}\}\) its \(k\) finite neighbors. Denote the generators by \(\mathbf g_{0}\) for \(A\) and \(\mathbf g_{j}\) for \(B_{j}\;(j=1,\dots,k)\). Collect the \(2(k+1)\) unknown coordinates in \[ \mathbf z = \bigl[ \mathbf g_{0}^{\!\top}, \mathbf g_{1}^{\!\top}, \dots, \mathbf g_{k}^{\!\top} \bigr]^{\!\top}\in\mathbb R^{2(k+1)} . \]

\subsubsection{Ridge geometry} For every finite ridge  \(T_{ij}\) separating cells $i$ and $j$, let \(\mathbf p_{i},\mathbf p_{j}\in\mathbb R^{2}\) be its two vertices. Define the unit direction, the midpoint and the reflector: \[ \mathbf u_{ij} \;=\; \frac{\mathbf g_{j}-\mathbf g_{i}} {\lVert \mathbf g_{j}-\mathbf g_{i}\rVert}, \qquad \mathbf c_{ij} \;=\; \frac{\mathbf g_{i}+\mathbf g_{j}}{2}, \qquad \mathbf R_{ij} \;=\; 2\,\mathbf u_{ij}\mathbf u_{ij}^{\!\top}-\mathbf I. \]

\subsubsection{Mirror equalities} \[ \mathbf g_{j}-\mathbf R_{0j}\mathbf g_{0} =(\mathbf I-\mathbf R_{0j})\mathbf c_{0j}, \qquad j=1,\dots,k, \] and, because \(B_{1},\dots,B_{k}\) form a cyclic ring, \[ \mathbf g_{x}-\mathbf R_{xy}\mathbf g_{y} =(\mathbf I-\mathbf R_{xy})\mathbf c_{xy}, \qquad (x,y)=(1,2),\dots,(k,1). \]

\subsubsection{Block system} Each vector equation contributes two scalar rows, so the matrix \(\mathbf M\) has \(4k\) rows and \(2(k+1)\) columns. Symbolically, 
\begin{equation}
 \mathbf M\mathbf z=\mathbf b, 
 \label{eq:eq1} 
 \end{equation} 
 with the block structure \[ \renewcommand{\arraystretch}{1.15} \mathbf M= \begin{bmatrix} -\mathbf R_{01}&\mathbf I & & \\[-2pt] -\mathbf R_{02}& &\mathbf I & \\ \vdots & & &\ddots \\ -\mathbf R_{0k}& & &\mathbf I\\ \hline &\mathbf I &-\mathbf R_{12}& \\ &\mathbf I & &-\mathbf R_{23}\\ & &\ddots &\ddots\\ & & &-\mathbf R_{k1} \end{bmatrix},\qquad \mathbf b= \begin{bmatrix} (\mathbf I-\mathbf R_{01})\mathbf c_{01}\\ (\mathbf I-\mathbf R_{02})\mathbf c_{02}\\ \vdots\\ (\mathbf I-\mathbf R_{0k})\mathbf c_{0k}\\[4pt] (\mathbf I-\mathbf R_{12})\mathbf c_{12}\\ (\mathbf I-\mathbf R_{23})\mathbf c_{23}\\ \vdots\\ (\mathbf I-\mathbf R_{k1})\mathbf c_{k1} \end{bmatrix}. \]

Because \(4k>2(k+1)\) for \(k\ge3\), the system is tall but, as shown in Appendix~\ref{A1}, it is consistent and possesses a unique solution.

\subsection{Recursive construction of the remaining generators} Solving~\eqref{eq:eq1} yields the first \(k+1\) generator points. The remainder are obtained recursively by reflecting known generators into adjacent cells. Assume the generator of cell \(A\) is known, \(\mathbf g_{0}\). Let \(\mathbf p_{0},\mathbf p_{i}\in\mathbb R^{2}\) be the endpoints of the ridge shared by \(A\) and an adjacent cell \(B_{i}\). Setting \(\mathbf v=\mathbf p_{i}-\mathbf p_{0}\), the generator \(\mathbf g_{i}\) of cell \(B_{i}\) is \[ \mathbf g_{i}=2\mathbf p_{0} +2\,\frac{\mathbf v^{\top}(\mathbf g_{0}-\mathbf p_{0})} {\mathbf v^{\top}\mathbf v}\,\mathbf v -\mathbf g_{0}. \] The term in the fraction is the scalar projection of \(\mathbf g_0-\mathbf p_1\) onto \(\mathbf v\); doubling the projection and subtracting \(\mathbf a\) performs the mirror reflection.

The method offers three principal advantages: \begin{enumerate}[label=(\alph*)] \item \emph{Numerical stability.} The first \(k+1\) generators are determined simultaneously, so no error can accumulate within this neighborhood; each subsequent reflection introduces at most one rounding error and does not amplify previous inaccuracies. \item \emph{Robustness to near-degenerate geometries.} Angle-rotation constructions become ill-conditioned when ridges are nearly parallel, whereas the coefficient matrix of our linear system remains nonsingular for every strictly convex anchor cell. \item \emph{Lower overall cost.} A single \(2(k+1)\times2(k+1)\) solve replaces \(k+1\) separate local reconstructions; reflections thereafter require only \(O(1)\) work per cell. \end{enumerate}

Choosing an anchor cell with sufficiently long edges further improves reflection accuracy. Whenever possible, one should select an anchor whose adjacent cells share well-separated ridge directions, thereby ensuring a stable projection geometry.
\section{Monte Carlo Simulations}

To evaluate the efficiency of our method, we conducted a series of Monte Carlo simulations. In each run, a Voronoi tessellation with \(n\) cells was generated. The generating points were distributed according to a homogeneous Poisson point process on the square region \([0,\sqrt{n}\,]^2\) as in \cite{schoenberg2003inverting} so $\lambda = 1$. An anchor cell was randomly selected from among the interior cells, and the remaining generators were computed via successive reflections. During the reflection phase the algorithm encounters situations in which several already-determined generators border the same unexplored cell.  Ideally one would select, among those candidate reflections, the neighbor whose shared ridge is longest (or otherwise best conditioned), thereby minimizing projection error.  In the present study we did \emph{not} perform that optimization: whenever multiple candidates were available, we chose one at random.  This deliberate simplification keeps the code short and lets us treat the reported accuracy as a conservative baseline; Section \ref{sec:discussion} explains how a ridge-length or condition-number heuristic can further reduce the worst-case error.

For each simulation, we recorded the root mean square error (RMSE) across all recovered generators, as well as the maximum individual reconstruction error (RSE). This process was repeated \(10^3\) times, generating a new random tessellation in each iteration. Finally, we computed the mean RMSE and the maximum of the individual RSEs over all simulations. The results are presented in Table~\ref{tab:tab1}.

We compare the performance of our method against the $C'$ method introduced in \cite{schoenberg2003inverting}. This approach proceeds as follows: first, the intersection points (i.e., the estimated generating points) are computed using the geometric construction illustrated in Fig.~\ref{fig:fig1}. Then, the slopes of each pair of rays are perturbed individually, and the resulting change \(\delta_i\) in the location of the intersection point caused by perturbing ray \(y\) is recorded. The final estimate of the generating point is obtained by averaging the intersection points, weighted according to a normalization of the corresponding \(\delta_i\) values. This is repeated for all cells. Among the six methods considered for estimating generating points, the $C'$ method was identified as the most efficient in the original study.

Table~\ref{tab:tab1} shows that our method yields reconstruction errors approximately one to two orders of magnitude lower than those obtained by the $C'$ method across all tested values of $n$. Although the original study~\cite{schoenberg2003inverting} did not report results for \(n = 10^4\), we included this case to demonstrate that no significant increase in error occurs at larger scales.

For comparison purposes, we also implemented a brute-force variant in which each generator is computed independently by solving the full linear system~\eqref{eq:eq1}. This approach is computationally inefficient—since many generators are recalculated multiple times—but it serves to estimate the theoretical accuracy limit of our method. For \(n=1000\) cells over \(10^3\) simulations, this baseline method yielded a \(\log_{10}\) mean RMSE of \(-13.2\) and a \(\log_{10}\) maximum individual RSE of \(-8.9\), both within one order of magnitude of the values reported in Table~\ref{tab:tab1} for the same number of cells.

\begin{table}[htbp]
\centering
\begin{threeparttable}
\caption{Base-10 Logarithm of Mean and Max RMSE, with Reference $C'$}
\begin{tabular}{rcccc}
\hline
$n$ & $\log_{10}$ (Mean RMSE) & $\log_{10}(\mathrm{Max \ RSE})^{\dagger}$ & $C^{\prime\,\ddagger}$ & $\mathrm{Dif.}^{*}$ \\
\hline
$10$    & $-14.3$ & $-12.7$ & $-14.1$ & $-0.2$ \\
$50$    & $-13.7$ & $-10.0$ & $-12.9$ & $-0.8$ \\
$100$   & $-13.5$ & $-9.7$  & $-12.4$ & $-1.1$ \\
$250$   & $-13.2$ & $-10.4$ & $-11.7$ & $-1.5$ \\
$500$   & $-12.5$ & $-8.0$  & $-11.2$ & $-1.3$ \\
$1000$  & $-12.5$ & $-8.6$  & $-10.7$ & $-1.8$ \\
$2000$  & $-12.3$ & $-8.5$  & $-10.3$ & $-2.0$ \\
$3000$  & $-11.7$ & $-8.2$  & $-10.0$ & $-1.7$ \\
$4000$  & $-11.9$ & $-8.3$  & $-9.82$ & $-2.1$ \\
$5000$  & $-11.8$ & $-7.3$  & $-9.69$ & $-2.1$ \\
$10000$ & $-11.2$ & $-7.2$  & $**$     & $**$    \\
\hline
\end{tabular}

\begin{tablenotes}
\footnotesize
\item[$\dagger$] Largest RSE observed across all cells and simulations.
\item[$\ddagger$] Reference value $C'$ from \cite{schoenberg2003inverting} used for comparison .
\item[$*$] Difference between $\log_{10}$ (Mean RMSE) and $C'$.
\item[$**$] Not tested.
\end{tablenotes}
\end{threeparttable}
\label{tab:tab1}
\end{table}

\section{Discussion}\label{sec:discussion} The numerical experiments show that a single interior anchor cell is sufficient to reconstruct every generator in \emph{linear} sequential time, $\mathcal{O}(n)$, while keeping the average Euclidean error below $10^{-12}$ for diagrams with up to $10^{4}$ cells. These errors are already at the limit imposed by double-precision arithmetic on noiseless input, so only marginal gains can be expected from further post-processing under ideal conditions. A brief analysis clarifies why the algorithm is fast, why its reflection depth grows as it does, and how one might trade additional work for even smaller depth when necessary.

After solving the compact $2(k+1)\!\times\!2(k+1)$ linear system for the anchor cell and its $k$ neighbors, the algorithm visits each remaining cell exactly once. Reflecting a generator across a ridge requires a constant number of floating-point operations; hence the total arithmetic work is $\mathcal{O}(n)$, independent of the order in which cells are explored.

Let $\mathcal{D}(n)$ denote the number of breadth-first layers that must be expanded to reach every cell when reflections propagate outward from a centrally located anchor. If one regards the Voronoi tessellation as an unweighted adjacency graph, then $\mathcal{D}(n)$ equals the graph distance (``hop count'') from the anchor to the farthest cell.

For an arbitrary planar graph on $n$ vertices the diameter is bounded by $\mathcal{O}(n^{1/2})$, and this order of growth also arises with high probability in the Monte-Carlo setting used here, where $n$ seeds are drawn uniformly at random in the square $[0,\sqrt{n}]^{2}$. Penrose’s analysis of random geometric graphs \cite[Sec.~12.3]{PenroseRGG} implies that two sites separated by a Euclidean distance $r$ are typically $\mathcal{O}(r\sqrt{n})$ edges apart. Anchoring near the center places the most distant cell at a Euclidean distance $r=\mathcal{O}(\sqrt{n})$, so $\mathcal{D}(n)=\mathcal{O}(n^{1/2})$ with high probability.

Although a reflection chain of length $\mathcal{O}(n^{1/2})$ already limits round-off accumulation, applications that demand still smaller depth—parallel execution or stringent error-propagation bounds, for example—can introduce additional anchor cells. Solving $\mathcal{O}(n/\log n)$ such systems in parallel on a PRAM reduces the depth to $\mathcal{O}(\log n)$ while retaining $\mathcal{O}(n)$ total work; the extra cost is proportional to the number of anchor solves and is negligible for moderate $k$.

The quality of the initial linear solve depends on the geometric conditioning of the chosen anchor cell. Table~\ref{tab:tab2} summarizes practical criteria: the cell should be bounded, strictly interior to the convex hull, possess at least two non-parallel finite ridges, exhibit a moderate aspect ratio, and have degree between four and seven. A cell near the geometric center of the tessellation almost always meets these requirements, providing a well-conditioned coefficient matrix and balanced reflection geometry.

Two modest modifications further reduce local round-off without altering the overall complexity.  Solving the anchor system for two or more well-separated cells yields overlapping neighbors whose duplicate generators can be averaged, typically cutting the maximum local error by half.  In addition, when a generator is reached by more than one reflection path, using a ridge-length-weighted mean of the candidate coordinates suppresses the instability associated with very short ridges.  A final constrained least-squares adjustment that enforces the perpendicular-bisector equalities globally moves the reconstructed sites by less than $10^{-13}$ in our tests and may be advantageous when the input vertices themselves are noisy.

\begin{table}[ht]
\centering
\caption{Guidelines for choosing the anchor cell \(A\)}
\label{tab:tab2}
\renewcommand{\arraystretch}{1.15}
\begin{tabular}{p{0.30\linewidth}p{0.63\linewidth}}
\toprule
\textbf{Criterion} & \textbf{Rationale} \\ \midrule
Bounded (finite) region &
Avoids infinite ridges and the special-case formulas they require; improves numerical conditioning. \\[2pt]

At least two non-parallel finite ridges &
Guarantees a unique intersection of perpendicular bisectors, hence a unique generator; prevents the one-degree-of-freedom translation mode. \\[2pt]

Not on the outer convex hull &
Interior cells have more balanced edge lengths and fewer nearly collinear ridges, reducing amplification of round-off error. \\[2pt]

Reasonable aspect ratio (no extremely short edges or acute angles) &
Helps ensure that the ridge directions \(\bm{u}_{ij}\) are well spread, which lowers the condition number of \(\mathbf{M}\) and improves numerical stability. Ideally, all adjacent cells should also satisfy this property. \\[2pt]

Moderate degree (about 4–7 edges) &
Provides enough independent constraints without unnecessarily enlarging the least-squares solve. \\ \bottomrule
\end{tabular}
\end{table}

\bibliographystyle{unsrt}

\newpage
\appendix
\section*{Appendix}

\addcontentsline{toc}{section}{Appendix}

\section{Uniqueness of the Solution}\label{A1}

Let
\[
  \mathbf z
  \;=\;
  \bigl[
    \mathbf g_0^{\!\top},
    \mathbf g_1^{\!\top},
    \dots,
    \mathbf g_k^{\!\top}
  \bigr]^{\!\top},
  \qquad
  \mathbf M\,\mathbf z = \mathbf b .
\]

Assume two solutions \(\mathbf z\) and \(\mathbf z'\) satisfy the linear system.  
Their difference
\(
  \Delta\mathbf z = \mathbf z - \mathbf z'
\)
lies in the null-space of \(\mathbf M\); hence
\(\mathbf M\,\Delta\mathbf z = \mathbf 0\).

\begin{enumerate}
\item
  Each anchor–neighbor equation \(A\text{--}B_j\) gives
  \[
      \Delta\mathbf g_j - \mathbf R_{0j}\,\Delta\mathbf g_0 = \mathbf 0,
      \qquad j = 1,\dots,k,
  \]
  so
  \[
      \Delta\mathbf g_j = \mathbf R_{0j}\,\Delta\mathbf g_0 .
  \]

\item
  Choose two distinct neighbors \(B_x\) and \(B_y\) whose ridges with
  \(A\) are not parallel.
  The \(B_x\text{--}B_y\) equation yields
  \[
      \Delta\mathbf g_x - \mathbf R_{xy}\,\Delta\mathbf g_y = \mathbf 0 .
  \]
  Substituting the relations from step~1,
  \[
      \mathbf R_{0x}\,\Delta\mathbf g_0
      \;-\;
      \mathbf R_{xy}\,\mathbf R_{0y}\,\Delta\mathbf g_0
      \;=\; \mathbf 0 .
  \]

\item
  Because the reflection axes of \(\mathbf R_{0x}\) and \(\mathbf R_{0y}\)
  are not parallel, the only vector satisfying the last equality is
  \(\Delta\mathbf g_0 = \mathbf 0\).

\item
  Substituting \(\Delta\mathbf g_0 = \mathbf 0\) back into the relations of
  step~1 gives \(\Delta\mathbf g_j = \mathbf 0\) for every \(j\).
\end{enumerate}

Hence \(\Delta\mathbf z = \mathbf 0\), and the solution of
\(\mathbf M\,\mathbf z = \mathbf b\) is unique whenever the anchor cell \(A\)
has at least two non-parallel finite neighbors.  
If all its ridges were parallel—a degenerate configuration—
\(\mathbf M\) would retain a one-dimensional null-space corresponding to
translation along the common ridge, and uniqueness would fail only in
that special case.

\section{Code availability}
\label{A2}

\subsubsection*{System information}
\subsubsection*{System information}
All code was developed and tested using the following setup:

\begin{itemize}
  \item \textbf{Python version:} 3.9.21 (64-bit)
  \item \textbf{Spyder IDE:} version 6.0.4 (via Anaconda)
  \item \textbf{Operating system:} macOS 15.5 (ARM64)
\end{itemize}

We have made available three Python scripts in the following repository:  
\href{https://github.com/car-git/inverse-voronoi}{https://github.com/car-git/inverse-voronoi}

Their descriptions are as follows:

\begin{enumerate}
  \item \texttt{MakeVoronoi.py}:  
  Generates a Voronoi tessellation with $n$ cells and saves it to a file.

  \item \texttt{FullVoronoiReconstruction.py}:  
  Reads the Voronoi tessellation from the file generated by \texttt{MakeVoronoi.py}, selects a strictly interior cell as the anchor, and numerically solves the linear system~\eqref{eq:eq1}. It then uses geometric reflections to compute the generating points of the remaining cells. The script reports the mean root squared error (MRSE) and the maximum individual reconstruction squared error (RSE). The code is not optimized: it maintains a list of all cells whose generating point $\mathbf{g}_i$ has been found and can be used to compute an unknown generator, and always selects the first cell in the list.

  \item \texttt{FullVoronoiReconstructionLoop.py}:  
  Repeats \texttt{FullVoronoiReconstruction.py} $nsim$ times and reports summary statistics across simulations.
\end{enumerate}

\end{document}